\documentclass{article}
\usepackage{amsmath,amsthm,amsfonts,amscd,eucal}

\numberwithin{equation}{section}

\def\cc{{\mathcal C}}

\def\cam{{\mathcal M}}

\def\cp{{\mathcal P}}

\def\cs{{\mathcal S}}

\def\ii{{\rm i\,}}

\def\bc{{\mathbb C}}

\def\bn{{\mathbb N}}
\def\br{{\mathbb R}}

\def\r{\rho}
\def\s{\sigma}

\def\f{\varphi}
\def\c{\chi}

\newtheorem{Thm}{Theorem}[section]
\newtheorem{Cor}[Thm]{Corollary}
\newtheorem{Prop}[Thm]{Proposition}

\theoremstyle{definition}
\newtheorem{Def}[Thm]{Definition}

\theoremstyle{remark}
\newtheorem{Rem}[Thm]{Remark} 
 
\begin{document}
\noindent
{\LARGE
\textbf{\textsf{Wigner-Yanase  information on quantum state space:}}}

\noindent
{\LARGE
\textbf{\textsf{the geometric approach}}}

\bigskip

Paolo Gibilisco\footnote[1]{Electronic mail: gibilisc@sci.unich.it}

\textit{Dipartimento di Scienze, Facolt\`a di Economia, Universit\`a
di Chieti-Pescara ``G. D'Annunzio'', Viale Pindaro 42, I--65127
Pescara, Italy.}

\bigskip

Tommaso Isola \footnote[2]{Electronic mail: isola@mat.uniroma2.it}

\textit{Dipartimento di Matematica, Universit\`a di Roma ``Tor 
Vergata'', Via della Ricerca Scientifica, I--00133 Roma, Italy.}

\markright{Wigner-Yanase information on quantum state space}


 \begin{abstract}
 In the search of appropriate riemannian metrics on quantum state
 space the concept of statistical monotonicity, or contraction under
 coarse graining, has been proposed by Chentsov.  The metrics with
 this property have been classified by Petz.  All the elements of this
 family of geometries can be seen as quantum analogues of Fisher
 information.  Although there exists a number of general theorems
 sheding light on this subject, many natural questions, also stemming
 from applications, are still open.  In this paper we discuss a
 particular member of the family, the Wigner-Yanase information. 
 Using a well-known approach that mimics the classical pull-back
 approach to Fisher information, we are able to give explicit formulae
 for the geodesic distance, the geodesic path, the sectional and
 scalar curvatures associated to Wigner-Yanase information.  Moreover
 we show that this is the only monotone metric for which such an
 approach is possible.
 \end{abstract}

  
 \section{Introduction} \label{sect.I}

 The notion of information proposed by Fisher is fundamental in
 probability and statistics for a number of reasons; here we mention
 only the Cramer-Rao inequality and the asymptotic behaviour of
 maximum lilkelihood estimators for exponential models (one can see
 \cite{Carlen} for unexpected features and applications of Fisher
 information).  In classical statistics Rao was the first to point out
 that Fisher information can be seen as a riemannian metric on the
 space of probability densities.  This point of view was nicely
 complemented by the results of Chentsov saying that (on the simplex
 of probability vectors) Fisher information is the unique riemannian
 metric contracting under Markov morphisms.  This can be rephrased in
 a more suggestive way.  Markov morphisms, or positive mappings, are
 the mathematical counterpart of the notion of noise.  Now suppose
 that we want to use a distance to distinguish different states
 (probability densities) in a statistically relevant way.  Then the
 effect of noise must be that of contracting the metric.  Chentsov
 theorem says therefore that in the classical case there is only one
 choice, the Fisher information (another argument producing Fisher
 information can be found in \cite{Wootters}).

 In the quantum case one deals with density operators instead of
 density vectors and completely positive mappings play the role of
 Markov morphisms.  As often happens in the quantum counterpart of a
 classical theory, instead of a uniqueness result, one has a
 classification theorem, due to Petz.  This result states that there
 is bijection between statistically monotone metrics on quantum state
 space and the operator monotone functions: we have therefore a rich
 ``garden" of candidates for the role of Fisher information in quantum
 physics.  Among the elements of this family of metrics one can find,
 in a certain sense, the most relevant riemannian metrics appeared in
 the literature \cite{Pe02}.

 Despite the existence of general results for the theory
 \cite{LeRu,Di00,GiIs2,Jen1,Jen2,Ha03,GiIs4} a number of open problems
 resists to investigation.  For example it does not exist yet a
 general formula for the geodesic path and the geodesic distance
 associated to an arbitrary monotone metric.  For the use of this kind
 of distances see for example \cite{NC}.  Because of the absence of a
 general formula, inequalities (giving bounds for the geodesic
 distance) have been proved \cite{Ruskai}.

 In this paper we discuss the Wigner-Yanase skew information.  To find
 the formulae for geodesic path and geodesic distance we mimic the
 classical approach to Fisher information via sphere geometry (one
 should note the importance of determining geodesic path in the study
 of the 2-Wasserstein metric \cite{CaGa}).  Indeed Wigner-Yanase
 information appears as the pull-back of the square root map
 \cite{GiIs3}.  Next we prove the formula for the scalar curvature. 
 One proof, due to J.Dittmann, uses the general formula \cite{Di00}
 and requires a long calculation.  The second one just uses the
 pull-back approach.  One should emphasize that, since the scalar
 curvature determines the asymptotic behaviour of the volume (for a
 riemannian metric) then it has also a statistical meaning in relation
 to the quantum analogue of Jeffrey's rule for determining prior
 probability distributions (see \cite{Pe02}).  Finally we prove, as a
 corollary of the results in \cite{HaPe,Ha03,GiIs4} that the
 Wigner-Yanase information is the only monotone metric that can be
 seen as a pull-back metric.

 The paper is organised as follows.  In section II we review the
 geometric approach to Fisher information.  In section III one finds
 an introduction to the general theory of statistical monotone
 metrics.  Section IV shows how the Wigner-Yanase information can be
 seen as a monotone riemannian metric.  In section V we show that the
 Wigner-Yanase geometry can be seen as the sphere geometry transposed
 on the space of density matrices; moreover we characterise it as the
 unique pull-back metric.  Section VI contains some comments on the
 main results and on some open problems.

 \section{Fisher information and its geometry} \label{sect.II}

 The classical definition of Fisher information for an indexed family
 of densities $p_{\theta}$ is given by the variance of the score.  In
 the case of a family indexed by only one parameter $\theta$ it is the
 number
 \begin{equation}
    I(\theta)=
    {\hbox {E}}_{\theta}\left[
        \left(\frac{\partial}{\partial \theta}  \log p_{\,\theta}\right)^2
    \right]
 \end{equation}
 assigned to the parameter $\theta$.  For $n$ parameters, say
 $\theta=\left(\theta^1,\dots, \theta^n\right)$, it is a matrix
 defined on the parameter manifold given by
 \begin{equation}
    I(\theta)_{i\,j}=
    {\hbox {E}}_{\theta}
        \left[
            \left(\frac{\partial}{\partial \theta^i}  \log p_{\,\theta}\right)
            \left(\frac{\partial}{\partial \theta^j}  \log p_{\,\theta}\right)
        \right] .
 \end{equation}
 Geometrically this means, that $I(\theta)$ is a symmetric bilinear
 form on the tangent spaces of the parameter manifold.  In a
 coordinate free language it reads as
 \begin{equation}\label{bilinearform}
     I(\theta)(U,V)= {\hbox {E}}_{\theta}\left[U\left(\log
     p_{\,\theta}\right)\, V \left(\log p_{\,\theta}\right)\right] ,
 \end{equation} 
 where $U$ and $V$ are vectors tangent to the parameter manifold and
 $U\left(\log p_{\,\theta}\right)$ is the derivative of $\log
 p_{\,\theta}$ along the direction $U$, that means $U\left(\log
 p_{\,\theta}\right)= {\frac{{\rm d}}{{\rm d} t}\log p_{\,\theta+t U}
 }_{|t=0}$.

 $I(\theta)$ is a measure for the statistical distinguishabilty of
 distribution parameters.  Under certain regularity conditions for
 $\theta \mapsto p_{\,\theta}$ the image of this mapping is a manifold
 of distributions.  This manifold is the actual object of interest in
 information geometry rather than the space of distribution parameters
 and formula (\ref{bilinearform}) defines a Riemannian metric $g$ on
 it (for a general reference see \cite{AmNa}).  Indeed, a vector $u$
 tangent to this manifold is of the form
 $$
 u={\frac{{\rm d}}{{\rm d} t}p_{\,\theta+t U} }_{\,|t=0}\,
 $$
 and the right hand side of (\ref{bilinearform}) now reads as
 \begin{equation}\label{Fishermetric}
 g(u,v):={\hbox {E}}_p\left[\frac{u}{p}\,\frac{v}{p}\right]
 \end{equation}
 defining the Fisher metric on the manifold of densities.  If the
 differential of $\theta \mapsto p_{\,\theta}$ is not injective, than
 there is some parameter redundance or ambiguity in the choice of $U$
 and $V$, and therefore the right hand side of (\ref{bilinearform})
 does not depend on this choice.

 We restrict now to $\cp_{n}\subset \br^{n}$, the simplex of strictly
 positive probability vectors, that is $\cp_{n} := \{ \r\in\br^{n} :
 \sum_{i=1}^{n} \r_{i} = 1, \r_{i}>0, i= 1,\ldots,n \}$.  An element
 $\r\in\cp_n$ is a density on the $n$-point set $\{1,\ldots,n\}$ with
 $\r(i)=\r_i$.  We regard an element $u$ of the tangent space ${\rm
 T}_\r\cp_{n}\equiv \{ u\in\br^{n} : \sum_{i=1}^{n} u_{i} = 0 \}$ as a
 function $u$ on $\{1,\ldots,n\}$ with $u(i)=u_i$.

\begin{Def}
The Fisher-Rao Riemannian metric on ${\rm T}_\r\cp_{n}$ is given by
\begin{equation}\label{FisherRao}
    \langle u,v {\rangle}^F_{\r} := \sum_{i=1}^{n} \frac{u_{i}v_{i}}{\r_{i}}
\end{equation}
for $u,\ v\in {\rm T}_{\r}\cp_{n}$.
\end{Def}

 To see the relation between this metric and the Fisher metric, let
 $u,v\in {\rm T}_\r\cp_{n}$.  We obtain from (\ref{Fishermetric})
 $$
 g(u,v)=\sum_{i=1}^{n} \frac{u(i)}{\r_i}\, \frac{v(i)}{\r_i}\r_{i}=
 \sum_{i=1}^{n} \frac{u_{i}v_{i}}{\r_{i}}
 $$
 in accordance with (\ref{FisherRao}).

 The following result is well known and is a very special case of a
 far more general situation (see \cite {Friedrich} for example).
\begin{Thm}
    The manifold $\cp_n$ equipped with the Fisher-Rao Riemannian
    metric $\langle\cdot,\cdot\rangle^F$ is isometric with an open
    subset of the sphere of radius 2 in $\br^{n}$.
\end{Thm}
\begin{proof}
    We consider the mapping $\f: \cp_{n}\to S^{n-1}_2\subset\br^n$,
    $$
    \f(\r):=2\left(\sqrt{\r_1},\dots,\sqrt{\r_n}\right)\, .
    $$
    Then
    $D_{\r}\f(u)=\left(\frac{u_1}{\sqrt{\r_1}},\dots,\frac{u_n}{\sqrt{\r_n}}\right)$
    and we get
    $$
    D_{\r}\f\left(\langle\cdot,\cdot\rangle^F \right)(u,v):= \langle
    D_{\r}\f (u),D_{\r}\f (v)\rangle^{\br^n} =\sum_{i=1}^{n} \frac{u_i
    \,v_i}{\r_{i}} =\langle u,v {\rangle}^F_{\r}\,.
    $$
    Hence the standard metric on the sphere of radius 2 is pulled back
    to the Fisher-Rao Riemannian metric.
\end{proof}

 This identification of $\cp_n$ with an open subset of a radius 2
 sphere allows for obtaining differential geometrical quantities of
 the Riemannian manifold
 $\left(\cp_n,\langle\cdot,\cdot\rangle^F\right)$.  From the very
 definition of geodesic distance, geodesic path and scalar curvature,
 one has for $S_r^{n-1}$, with $P_1, P_2\in S_r^{n-1}$,
 \newline
 1) geodesic distance
 $$
 d(P_1,P_2) = r \cdot \hbox{arcos}\left(\frac{\left\langle P_1,P_2
 \right\rangle}{r^2}\right)
 $$
 \newline
2) geodesic path connecting $P_1$ and $P_2$ :
 $$
 \gamma^{P_1,P_2}(t)= r\frac {(1-t) P_1+t P_2} {||(1-t) P_1+t P_2 ||}
 $$
 (of course, $t$ is not the arc length parameter);
\newline
3) scalar curvature
 $$
 \hbox{Scal}(v)= \frac{1}{r^2}(n-1)(n-2)
 $$
 because  $S^{n-1}_r$ has constant sectional curvature equal to $\frac{1}{r^2}$.

 Let us denote by $d_F,\gamma_F, \hbox{Scal}_F$ respectively the
 corresponding quantities for the Fisher information.  The above
 considerations give, for $\rho,\s\in \cp_n$,
\newline
1) Bhattacharya distance
$$
d_{\hbox{F}}({\rho},\s) = 2\hbox{arccos}\left(\sum_i {\rho}^{1/2}_i
{\s}^{1/2}_i\right)
$$
\newline
2) geodesic path connecting $\rho$ and $\s$:
$$
\gamma_{\hbox{F}}^{\rho,\s}(t)= 2\frac {((1-t)
\sqrt{\rho}+t\sqrt{\s})^2} {\sum_i ((1-t)
\sqrt{\rho_i}+t\sqrt{\s_i})^2}
$$
\newline
3) scalar curvature
$$
\hbox{Scal}_{\hbox{F}} (\rho)=\frac{1}{4}(n-1)(n-2) \qquad \qquad
\forall \rho \in \cp_{n}.
$$

The Levi-Civita connection associated to Fisher metric can be
decomposed using the geometry of mixture and exponential models.  The
rest of the section explains how.

\begin{Def}\label{dualistic}
A dualistic structure on a manifold $\cam$ is a triple $(\langle \cdot,\cdot \rangle, \nabla, \widetilde{\nabla})$ where $\langle \cdot,\cdot \rangle$ is a riemannian metric on $\cam$ and $\nabla, \widetilde{\nabla}$ are affine connections on $\cam$ such that
$$
X\langle Y,Z \rangle=\langle \nabla_XY,Z \rangle + \langle
Y,\widetilde{\nabla}_XZ \rangle
$$
where $X,Y,Z$ are vector fields.  If
$U^{\nabla},U^{\widetilde{\nabla}}$ are the parallel transport
associated to $\nabla, \widetilde{\nabla}$ then the above equation is
equivalent to
$$
\langle U^{\nabla}(u),U^{\widetilde{\nabla}}(v) \rangle = \langle u,v
\rangle.
$$

\end{Def}

A divergence on a manifold is a smooth nonnegative function $ D: \cam
\times \cam \to \br$ such that $D(\r,\s)=0$ iff $\r=\s$.  To each
divergence $D$ one may associate a dualistic structure $(\langle
\cdot,\cdot \rangle, \nabla, \widetilde{\nabla})$ (see
\cite{AmNa,Eguchi}).

Let $\nabla^2$ be the Levi-Civita connection of Fisher information. 
The Kullback-Leibler relative entropy $K(\rho,\s)=\sum_i \rho_i (\log
\rho_i - \log \s_i)$ gives a dualistic structure
$(\langle\cdot,\cdot\rangle^F, \nabla^m,\nabla^e)$ such that
$$
\nabla^2=\frac{1}{2}(\nabla^m+\nabla^e)
$$
where $\nabla^m,\nabla^e$ are the mixture and exponential connections. 
These connections are torsion free and flat: once the representation
by scores is used for the tangent spaces, the associated parallel
transports are given by
$$
U^m_{\rho \s}: T_{\r}{\cal P} \to T_{\s}{\cal P} \qquad \qquad
U^m_{\rho \s}(u)= \frac{\rho}{\s}u
$$
$$
U^e_{\rho \s}: T_{\r}{\cal P} \to T_{\s}{\cal P} \qquad \qquad
U^e_{\rho \s}(u)= u-\hbox{E}_\s(u).  
$$
The geodesics of $\nabla^m,\nabla^e$ are, respectively, the mixture
and exponential models.

\section {Metric contraction under coarse graining} \label{sect.III}

In the commutative case a Markov morphism (or stochastic map) is a
positive operator $T: \br^n \to \br^k$.  In the noncommutative case a
stochastic map is a completely positive and trace preserving operator
$T: M_n \to M_k$ where $M_n$ denotes the space of $n$ by $n$ complex
matrices.  We shall denote by ${\cal D}_n$ the manifold of strictly
positive elements of $M_n$ and by ${\cal D}^1_n \subset {\cal D}_n$
the submanifold of density matrices.

In the commutative case a monotone metric is a family of riemannian metrics
$g=\{g^n\}$ on $\{\cp_n\}$, $n \in \bn$ such that
$$
g^m_{T(\rho)}(TX,TX) \leq g^n_{\rho}(X,X)
$$
holds for every stochastic mapping $T:\br^n \to \br^m$ and all $\rho
\in \cp_n$ and $X \in T_\rho \cp_n$.

In perfect analogy, a monotone metric in the noncommutative case is a
family of Riemannian metrics $g=\{g^n\}$ on $\{{\cal D}^1_n\}$, $n \in
\bn$ such that
$$
g^m_{T(\rho)}(TX,TX) \leq g^n_{\rho}(X,X)
$$
holds for every stochastic mapping $T:M_n \to M_m$ and all $\rho \in
{\cal D}^1_n$ and $X \in T_\rho {\cal D}^1_n$.

 Let us recall that a function $f:(0,\infty)\to \br$ is called
 operator monotone if for any $n\in \bn$, any $A$, $B\in M_n$ such
 that $0\leq A\leq B$, the inequalities $0\leq f(A)\leq f(B)$ hold. 
 An operator monotone function is said symmetric if $f(x):=xf(x^{-1})$
 and normalized if $f(1)=1$.  In what follows by operator monotone we
 mean normalised symmetric operator monotone.  With each operator
 monotone function $f$ one associates also the so-called
 Chentsov--Morotzova function
 $$
 c_f(x,y):=\frac{1}{yf(\frac{x}{y})}\qquad\hbox{for}\qquad
 x,y>0.
 $$
 Define $L_{\r}(A) := \r A$, and $R_{\r}(A) := A\r$.  Since $L_{\r},
 R_{\r}$ commute we may define $c(L_{\r}, R_{\r})$.  Now we can state
 the fundamental theorems about monotone metrics (uniqueness and
 classification are up to scalars).

\begin{Thm} {\rm \cite{Chentsov}}
    There exists a unique monotone metric on $\cp_n$ given by the
    Fisher information.
\end{Thm}

 \begin{Thm}{\rm \cite{Pe96}}
     There exists a bijective correspondence between monotone metrics
     on ${\cal D}_n^1$ and operator monotone functions given by the
     formula
     $$
     \langle A,B {\rangle}_{\rho,f}:=Tr(A\cdot c_f(L_\rho,R_\rho)(B)).
     $$
 \end{Thm}

 The tangent space to ${\cal D}^1_{n}$ at $\r$ is given by
 $T_{\r}{\cal D}^1_{n} \equiv \{ A\in M_{n} : A=A^{*}, Tr(A)=0\}$, and
 can be decomposed as $T_{\r}{\cal D}^1_{n} = (T_{\r}{\cal
 D}^1_{n})^{c} \oplus (T_{\r}{\cal D}^1_{n})^{o}$, where $(T_{\r}{\cal
 D}^1_{n})^{c}:= \{ A\in T_{\r}{\cal D}^1_{n} : [A,\r] = 0\}$, and
 $(T_{\r}{\cal D}^1_{n})^{o}$ is the orthogonal complement of
 $(T_{\r}{\cal D}^1_{n})^{c}$, with respect to the Hilbert-Schmidt
 scalar product $\langle A,B \rangle := Tr(A^{*}B)$.  Each
 statistically monotone metric has a unique expression (up to a
 constant) given by $Tr(\r^{-1}A^{2})$, for $A\in (T_{\r}{\cal
 D}^1_{n})^{c}$.  The following result will be used in Section
 \ref{sect.V}.

 \begin{Prop}\label{derivative}(See {\rm \cite{Bhatia}}). 
     Let $A \in T_{\r}{\cal D}^1_{n}$ be decomposed as $A=A^c+i[\r,U]$
     where $A^c \in (T_{\r}{\cal D}^1_{n})^{c}$ and $i[\r,U] \in
     (T_{\r}{\cal D}^1_{n})^{o}$.  Suppose $\f \in {\cal
     C}^1(0,+\infty)$.  Then
     $$
     (D_{\r}\f)(A)=\f'(\r)A^c+i[\f(\r),U].
     $$
 \end{Prop}

 As proved by Lesniewski and Ruskai each monotone metric is the
 hessian of a suitable relative entropy; to state this result more
 precisely, we introduce some notation.  In what follows $g$ is an
 operator convex function defined on $(0,+\infty)$ and such that
 $g(1)=0$.  The formula
     $$
     f(x)\equiv f_g(x):= \frac{(x-1)^2}{g(x)+xg(x^{-1})}
     $$
 associates a normalised, symmetric operator monotone function $f=f_g$
 to each $g$.  We denote by $\Delta_{\s,\rho}=L_{\s}R^{-1}_{\rho}$ the
 relative modular operator.  The relative $g$--entropy of $\rho$ and
 $\s$ is defined as
 $$
     H_g(\rho,\s) :=
     Tr(\rho^{\frac{1}{2}}g(\Delta_{\s,\rho})(\rho^{\frac{1}{2}})).
 $$
 $H_g$ is a divergence on ${\cal D}_n$ in the sense of
 \cite{Eguchi,AmNa}.  If $\rho,\s$ are diagonal $H_g$, reduces to the
 commutative relative $g$--entropy (see \cite{Csiszar}).

\begin{Thm} {\rm \cite{LeRu}} 
    Let $g$ be operator convex, $g(1)=0$, $f=f_g$ and $\r \in {\cal
    D}_n$.  Then
    $$
    -\frac{\partial}{\partial t}\frac{\partial}{\partial
    s}H_g(\r+tA,\r+sB){\Bigm |}_{t=s=0}
    =Tr(A\cdot c_f(L_\rho,R_\rho)(B)).
    $$
\end{Thm}

\medskip

 To state the general formula for the scalar curvature of a monotone
 metric we need some auxiliary functions.  In what follows $c', (\log
 c)'$ denote derivatives with respect to the first variable, and
 $c=c_f$.

\begin{eqnarray}
 h_1(x,y,z)&:=&\frac{c(x,y)-z\,
 c(x,z)\,c(y,z)}{(x-z)(y-z)c(x,z)c(y,z)}\,,\nonumber\\
 h_2(x,y,z)&:=&\frac{\left(c(x,z)-c(y,z)\right)^2}{(x-y)^2c(x,y)c(x,z)c(y,z)}\,,\nonumber\\
 h_3(x,y,z)&:=&z\,\frac{(\ln c)'(z,x)-(\ln c)'(z,y)}{x-y}\,,\nonumber\\
 h_4(x,y,z)&:=&z\,(\ln c)'(z,x)\;(\ln c)'(z,y)\,,\nonumber\\
 h &:=& h_1-\frac{1}{2}\,h_2+2h_3-h_4\,\label{function_h}.
\end{eqnarray}
The functions $h_i$ have no essential singularities if arguments coincide.

Note that $\langle A,B {\rangle}^f_{\rho}:=Tr(A\cdot
c_f(L_\rho,R_\rho)(B))$ defines a riemannian metric also over ${\cal
D}_n$ (${\cal D}^1_n$ is a submanifold of codimension 1).  Let
$\hbox{Scal}_f(\rho)$ be the scalar curvature of $({\cal D}_n,\langle
\cdot,\cdot {\rangle}^f_{\rho})$ at $\rho$ and $\hbox{Scal}^1_f(\rho)$
be the scalar curvature of $({\cal D}^1_n,\langle \cdot,\cdot
{\rangle}^f_{\rho})$.

\begin{Thm}\label{sc}{\rm \cite{Di00}}
Let $\sigma(\rho)$ be the spectrum of $\rho$. Then
\begin{equation}
    \hbox{Scal}_f(\rho)= \sum_{x,y,z \in \sigma(\rho)} h(x,y,z)-
    \sum_{x \in \sigma(\rho)} h(x,x,x)\label{scalcurv}
\end{equation}
$$
\hbox{Scal}^1_f(\rho)=\hbox{Scal}_f(\rho)+ \frac{1}{4}(n^2-1)(n^2-2).
$$
\end{Thm}

\section{Wigner-Yanase information as a riemannian metric} \label{sect.IV}

\bigskip

 Let $\r\in {\cal D}^1_{n}$ be a density matrix and let $A$ be a self
 adjoint matrix.  The Wigner-Yanase information (or skew information,
 information content relative to $A$) was defined as
 $$
 I(\r,A) := -\hbox{Tr}([\r^{1/2},A]^{2})
 $$
 where $[\cdot,\cdot]$ denotes the commutator(see \cite{WiYa}). 
 Consider now $g(x):=g_{\hbox{wy}}(x):=4(1-\sqrt{x})$.  In this case
 $$
 H_g(\rho,\s)=4(1- \hbox{Tr}(\rho^{\frac{1}{2}}\s ^{\frac{1}{2}})).
 $$
 The associated operator monotone and Chentsov-Morotzova functions are
 $$
 f_{\hbox{wy}}(x) := \frac {1}{4}(\sqrt{x}+1)^{2} \qquad \qquad
 c_{\hbox{wy}}(x,y) := \frac{1}{yf_{\hbox{wy}}(\frac{x}{y})} =
 \frac{4}{(\sqrt{x}+\sqrt{y})^{2}}
 $$

 Let us consider the monotone metric
 $$
  \langle A,B {\rangle}^{\hbox{wy}}_{\rho} :=
  {\hbox {Tr}}(A\,c_{\hbox{wy}}(L_{\r},R_{\r})(B)).
 $$
 A typical element of $(T_{\r}D_{n})^{o}$ has the form $\ii[\r,A]$,
 where $A$ is self-adjoint.  We have
 \begin{align*}
    \langle \ii[\r,A],\ii[\r,A]  {\rangle}^{\hbox{wy}}_{\rho} &
    =
    \hbox{Tr\,}\left(\ii[\r,A]4 (L_\r^{1/2}+R_\r^{1/2})^{-2}(\ii[\r,A])\right) \\
    & =
     - 4\,\hbox{Tr\,}\left( (L_\r^{1/2}+R_\r^{1/2})^{-1}([\r,A])\,\,\,
    (L_\r^{1/2}+R_\r^{1/2})^{-1}([\r,A]) \right) \\
    & =
    - 4\,\hbox{Tr\,}\left( (L_\r^{1/2}+R_\r^{1/2})^{-1}\circ(L_\r-R_\r)(A)\,\,\,
    (L_\r^{1/2}+R_\r^{1/2})^{-1}\circ(L_\r-R_\r)(A) \right) \\
    & =
      - 4\,\hbox{Tr\,}\left( (L_\r^{1/2}-R_\r^{1/2})(A)\,\,\,
    (L_\r^{1/2}-R_\r^{1/2})(A) \right) \\
    & =
     -4\,\hbox{Tr}\left([\r^{1/2},A]^{2}\right) = 4I(\r,A)
 \end{align*}
 and this explains why the monotone metric associated with the
 function $\frac {1}{4}(\sqrt{x}+1)^{2}$ is called the Wigner-Yanase
 monotone metric.

\section{The main result} \label{sect.V}

\bigskip

First of all we calculate the scalar curvature of Wigner-Yanase
information using Theorem \ref{sc}.  If $ f_{\hbox{wy}}(x) := \frac
{1}{4}(\sqrt{x}+1)^{2}$ we write $\hbox{Scal}^1_{\hbox{wy}}$ for
$\hbox{Scal}^1_f$.

\begin{Thm}\label{Jochen}
 $$
\hbox{Scal}^1_{\hbox{wy}}(\rho)= \frac{1}{4}(n^2-1)(n^2-2).
$$
\end{Thm}
\begin{proof}
Let us calculate the auxiliary functions for $c_{\hbox{wy}}(x,y) :=
4(\sqrt{x}+\sqrt{y})^{-2}$.  We get
\begin{eqnarray*}
h_1(x,y,z)&=& \frac{{\sqrt{x}}\,{\sqrt{y}} + 3\,{\sqrt{x}}\,{\sqrt{z}} +
     3\,{\sqrt{y}}\,{\sqrt{z}} + z}{4\,{{\left( {\sqrt{x}} + {\sqrt{y}} \right) }^2}\,
     \left( {\sqrt{x}} + {\sqrt{z}} \right) \,
     \left( {\sqrt{y}} + {\sqrt{z}} \right) }  \\
h_2(x,y,z)&=&\frac{{{\left( {\sqrt{x}} + {\sqrt{y}} + 2\,{\sqrt{z}} \right) }^2}}
   {4\,{{\left( {\sqrt{x}} + {\sqrt{z}} \right) }^2}\,
     {{\left( {\sqrt{y}} + {\sqrt{z}} \right) }^2}}\\
h_3(x,y,z)&=&\frac{{\sqrt{z}}}{\left( {\sqrt{x}} + {\sqrt{y}} \right) \,
     \left( {\sqrt{x}} + {\sqrt{z}} \right) \,
     \left( {\sqrt{y}} + {\sqrt{z}} \right) }\\
h_4(x,y,z)&=& \frac{1}{\left( {\sqrt{x}} + {\sqrt{z}} \right) \,
     \left( {\sqrt{y}} + {\sqrt{z}} \right). }
\end{eqnarray*}

 Now one can verify by calculation that the symmetrization of
 $h_1-\frac{1}{2}\,h_2$ and the symmetrization of $2\,h_3-h_4$ vanish. 
 Hence, by (\ref{function_h}), the symmetrization of $h$ vanishes,
 too.  Since we sum up in formula (\ref{scalcurv}) over all triples of
 eigenvalues we may replace $h$ with its symmetrization without
 changing the summation result.  Therefore
 $$
 \hbox{Scal}_{\hbox{wy}}(\r)=0\,,\qquad\qquad
 \hbox{Scal}^1_{\hbox{wy}}(\rho)= \frac{1}{4}(n^2-1)(n^2-2) \qquad
 \qquad \forall \r \in {\cal D}_n^1.
 $$
\end{proof}

In what follows we use the pull-back approach to derive (and explain)
the above formula in a direct way.  Furthermore we deduce the geodesic
distance and geodesic equation.

Let us denote by ${\cal S}$ the manifold $\left\{ A\in M_n:\,{\rm
Tr}\, A\,A^\ast=2, A=A^\ast \right\}$.  Clearly, since $\cal S$ is the
intersection of the radius 2 sphere in $\bc^{n\times n}$ and the
subspace of Hermitian matrices, it is isometric with a radius 2 sphere
$S^{n^2-1}_2$.

Now, let $\f : {\cal D}^1_{n} \to {\cal S}\subset \bc^{n\times n}$,
$\f(\r):=2 \sqrt{\r}$.  Then we have the following result (see
\cite{HaPe,GiIs3,Jen1,Grasselli}).

\begin{Thm}
 The pull-back by the map $\f$ of the natural metric on ${\cal
 S}\equiv S^{n^2-1}_2$ coincides with the Wigner-Yanase monotone
 metric.
\end{Thm}
\begin{proof}
 Let $A$ and $B$ be vectors tangent to ${\cal D}^1_n$ at $\r$. 
 Because $\f(\r)\,\f(\r)=4\,\r$ we get from the Leibniz rule
 $D_{\r}\f(A)\sqrt{\r}\,+\sqrt{\r}\,D_{\r}\f(A)=2A$ Thus, the
 differential of $\f$ at the point $\r$ is given by
 $$
 D_{\r}\f(A)= 2\left(L_{\r}^{1/2}+ R_{\r}^{1/2}\right)^{-1}\!(A)\,.
 $$
 Therefore the pull-back of the real part of the Hilbert-Schmidt metric yields
 \begin{align*}
     D_{\r}\f({\rm Re\,}\langle\cdot,\cdot\rangle)(A,B) & ={\rm
     Re\,}\langle D_{\r}\f (A),D_{\r}\f(B)\rangle \\
     & =4\,{\rm Re\,}\langle(L_{\r}^{1/2}+ R_{\r}^{1/2})^{-1}(A),
     (L_{\r}^{1/2}+ R_{\r}^{1/2})^{-1}\!(B)\rangle \\
     & =4\,\langle A,(L_{\r}^{1/2}+ R_{\r}^{1/2})^{-2}(B)\rangle  \\
     & = 4\,\hbox {Tr}\,A\,(L_{\r}^{1/2}+ R_{\r}^{1/2})^{-2}(B) \\
     & = \hbox {Tr}\,A \,c_{wy} (L_{\r},R_{\r})(B) = \langle A,B
     {\rangle}^{\hbox{wy}}_{\rho}
\end{align*}
    which was to be proved.
 \end{proof}

From this result one can deduce the following

\begin{Thm}
    For the geodesic distance, geodesic path and the scalar curvature
    of Wigner-Yanase information the following formulae hold

1) geodesic distance

\begin{equation}\label{eq:geodist}
d_{\hbox{wy}}({\rho},\s)=2\hbox{arccos}({\hbox{Tr}}( {\rho}^{1/2} {\s}^{1/2}))
\end{equation}

2)  geodesic path

\begin{equation}\label{eq:geopath}
{\gamma}_{\hbox{wy}}^{\rho,\s}(t)=
2\frac {((1-t) \sqrt{\rho}+t\sqrt{\s})^2} {{\hbox{Tr}}(((1-t) \sqrt{\rho}+t\sqrt{\s})^2)}
\end{equation}

3) scalar curvature

 \begin{equation}\label{eq:sc}
\hbox{Scal}^1_{\hbox{wy}}(\rho)= \frac{1}{4}(n^2-1)(n^2-2).
\end{equation}
 \end{Thm}

\begin{proof}
The formulae are immediate consequences of the preceding theorem and
of sphere geometry.  Indeed by the pull-back argument the
Wigner-Yanase metric looks locally like a sphere of radius 2 of
dimension $(n^2-1)$.  But for a sphere of this kind the sectional
curvatures are all equal to $\frac{1}{4}$ and therefore the scalar
curvature is given by $\frac{1}{4}(n^2-1)(n^2-2)$.
\end{proof}

One may ask if other monotone metrics are the pull-back of some
function $\f$ different from the square root.  The rest of the section
answers this question.

\begin{Def}
A monotone metric $\langle \cdot,\cdot \rangle_{\r,f}$ is a pull-back
metric if there exists a manifold $\cs \subset M_n$ and a function $\f
\in\cc^{1}(0,+\infty)$ such that the pull-back metric of $\f:{\cal
D}^1_n \to \cs \subset M_n$ coincides with $\langle \cdot,\cdot
\rangle_{\r,f}$.
\end{Def}

\begin{Prop}
Let $\langle \cdot,\cdot \rangle_{\r,f}$ be a monotone metric, let
$c=c_f$ be the associated $CM$-function and let $\f
\in\cc^{1}(0,+\infty)$.  We have that $\langle \cdot,\cdot
\rangle_{\r,f}$ is a pull-back metric by $\f$ if and only if
\begin{equation}\label{eq:pullback}
	 \left(\frac{\f(x)-\f(y)}{x-y}\right)^2= c(x,y).
\end{equation}
\end{Prop}
\begin{proof}
Apply the formula (\ref{derivative}) to tangent vectors in
$(T_{\r}{\cal D}^1_{n})^{o}$.
\end{proof}

\begin{Def}
Let $\f,\c\in\cc^{1}(0,+\infty)$.  We say that $(\f,\c)$ is a dual
pair if there exist an operator monotone $f$ such that
\begin{equation}\label{eq:dualpair}
\frac{\f(x)-\f(y)}{x-y}\cdot \frac{\c(x)-\c(y)}{x-y}= c(x,y).
\end{equation}
where $c=c_f$ is the $CM$-function associated to $f$.
\end{Def}

In such a case we say that $f$ (or $c_f$) is a dual function.  If
$(\f,\f)$ is a dual pair with respect to $f$ (or $c_f$) we say that
$f$ (or $c_f$) is self-dual.  Obviously one has

\begin{Prop}
To say that $\langle \cdot,\cdot \rangle_{\r,f}$ is a pull-back metric
by $\f$ it is equivalent to say that $f$ (or $c_f$) is self-dual with
respect to $\f$.
\end{Prop}

\begin{Def}
     Two dual pairs $(\f,\c), (\tilde {\f}, \tilde {\c})$ are
     equivalent if there exist constants $A_1,A_2,B_1,B_2$ such that
     $A_1A_2=1$
     \begin{align*}
	 \tilde {\f} & = A_1\f+B_1 \\
	 \tilde {\c} & = A_2\c+B_2.
     \end{align*}
 \end{Def}
 
 Obviously equivalent pairs define the same $CM$-function.  In what
 follows we consider dual pairs up to this equivalence relation with
 the traditional abuse of language.
 We are ready to state the fundamental result of the theory that classifies dual pairs.

 \begin{Thm} \label{main} {\rm (\cite{HaPe,Ha03,GiIs4})}
     Let $\f,\c\in\cc^{1}(0,+\infty)$.  Then $(\f,\c)$ is a dual pair
     if and only if one of the following two possibilities hold
     $$
     (\f(x),\c(x))=(\frac{x^p}{p},\frac{x^{1-p}}{1-p}) \qquad 
     p\in[-1,2]\setminus\{0,1\}
     $$
     $$
     (\f(x),\c(x))=(x,\hbox{\rm{log}}(x)).
     $$
 \end{Thm}

 \begin{Cor}
     The function $f(x)=\frac {1}{4}(\sqrt{x}+1)^{2}$ is the only
     self-dual operator monotone function, that is: the Wigner-Yanase
     metric is the only pull-back metric among statistically monotone
     metrics.
\end{Cor}

\section{Conclusions} 

\begin{Rem}
 Note that the formula (\ref{eq:geodist}) implies
 $d^{\hbox{wy}}({\rho},\s) \leq 2\pi$.  An analogous inequality holds
 for the Bures metric (see \cite{Di95}, p.311).  It seems that in the
 literature there are no other explicit formulas for the geodesics
 distance.  For example it is known that the formula
 \begin{equation}\label{eq:Bures}
     d_{\hbox{b}}(\rho,\s)= \sqrt{2-2 \hbox{Tr}(\rho^{ \frac{1}{2}} \s
     \rho^{\frac{1}{2}})^\frac{1}{2}}
 \end{equation}
 defines a metric on the state space whose infinitesimal counterpart
 (say the hessian) is the $SLD$-metric (that is
 $f(x)=\frac{1}{2}(1+x)$).  But this does not imply that Equation
 (\ref{eq:Bures}) is the geodesic distance of the $SLD$-metric.
\end{Rem}

\begin{Rem}
 In general it is difficult to give explicit formulae for geodesic
 paths of monotone metrics.  In the case of the Bures metric these
 geodesics can be given because they are projections of large circles
 on a sphere in the purifying space (see \cite{Di95} p.311 and
 \cite{DiUh}\cite{BrCa}).  For a discussion of geodesics for
 $\alpha$-connections see \cite{Jen1,Jen2}.
\end{Rem}

\begin{Rem}
 A classical theorem classifies the spaces of costant curvature
 \cite{KoNo}.  It is not known at the moment if there are other
 monotone metrics of costant sectional and scalar curvature.
\end{Rem}

\begin{Rem}
 We have seen in the commutative case that for the Levi-Civita
 connection of the pull-back of the square root it is available the
 decomposition
 $$
 \nabla^2=\frac{1}{2}(\nabla^m+\nabla^e).
 $$
 In the non-commutative case an analogous decomposition for the
 pull-back of the square root no longer holds.  Indeed, on one hand,
 the use of Umegaki relative entropy $H(\rho,\s)=\hbox{Tr}(\rho (\log
 \rho-\log \s))$ produces a similar decomposition, but for the
 Bogoliubov-Kubo-Mori metric \cite{Nagaoka,AmNa,GrSt2}.  On the other
 hand, if one uses
 $H_{\hbox{wy}}(\rho,\s)=4(1-\hbox{Tr}(\rho^{1/2}\s^{1/2}))$ as a
 divergence on $ {\cal D}^1_n$ and constructs the associated dualistic
 structure $(\langle \cdot,\cdot
 \rangle^{H_{\hbox{wy}}},\nabla^{H_{\hbox{wy}}},
 \nabla^{H_{\hbox{wy}}})$ (again following the lines of
 \cite{Eguchi,AmNa}), then the construction is trivial, namely the
 dual connections both coincide with the Levi-Civita connection of the
 Wigner-Yanase information.  This is easily seen on ${\cal P}_n$ where
 $H_g(\rho,\s)$ reduces to Csiszar relative $g$-entropy: it is known
 that such an entropy induces the $\alpha$-geometry where $\alpha$ is
 given by the formula $\alpha= 3+2g'''(1)/g''(1)$ (see \cite{AmNa}
 p.57).  For $g=4(1-\sqrt{x})$ this gives $\alpha=0$ that is the
 Fisher information case (see also \cite{Grasselli}).
\end{Rem}

   $$ $$
\noindent
 {\large \textbf{\textsf{ACKNOWLEDGEMENTS}}}

\bigskip This research has been supported by the italian MIUR program
"Quantum Probability and Infinite Dimensional Analysis" 2001-2002.  It
is a pleasure to thank J. Dittmann for a number of valuable
conversations on this subject and for permitting us to reproduce in
this paper the proof of Theorem \ref{Jochen}.



\end{document}